
\documentstyle{amsppt}

\nologo
\magnification=\magstep1
\hsize=145mm
\vsize=220mm
\hcorrection{-7mm}
\vcorrection{-10mm}
\baselineskip=18pt
 \abovedisplayskip=4pt
 \belowdisplayskip=4pt
\parskip=3pt
\parindent=8mm


\define\dnl{\newline\newline}

\define\nl{\newline}

\define\dg{{\roman{deg}}}
\define\dm{{\roman{dim}}}

\define\rk{{\roman{rank}}}

\define\SA{\Cal A}
\define\SE{\Cal E}
\define\SN{\Cal N}
\define\SO{\Cal O}
\define\BP{{\Bbb P}}

\define\Conj{{\bf Conjecture. }}
\define\Cor{{\bf Corollary. }}

\define\Fact{{\it Fact. }}

\define\Th{{\bf Theorem. }}

\topmatter
\author Takao FUJITA \endauthor
\address {Takao FUJITA\newline
Department of Mathematics \newline
Tokyo Institute of Technology \newline
Oh-okayama, Meguro, Tokyo \newline
152 Japan \newline
e-mail:fujita\@math.titech.ac.jp}
\endaddress
\title A note on scrolls of smallest embedded codimension
\endtitle
\endtopmatter

\document

\noindent {\bf \S0. Introduction}

The situation considered in this note is as follows.
Let $M$ be an algebraic submanifold of $\BP^N$ with $n=\dm M$.
$M$ is said to be a scroll over $S$ if there is a surjective morphism $\pi:
M\to S$ such that every fiber $F_x=\pi^{-1}(x)$ over $x\in S$ is a linear
subspace in $\BP^N$ of dimension $r-1=n-s$, where $s=\dm S$.
This is equivalent to saying that $M\cong\BP_S(\SE)$ for some vector bundle
$\SE$ of rank $r=n-s+1$ and the tautological bundle $H(\SE)$ is the hyperplane
section bundle of $M$.

We have $b_i(M)=b_i(\BP^N)$ for $i\le 2n-N$ by Barth-Lefschetz Theorem, hence
$N\ge 2n-1$ for scrolls, since otherwise $1+b_2(S)=b_2(M)=b_2(\BP^N)=1$,
contradiction.
Thus we want to study scrolls such that $N=2n-1$.

The case $s=1$ was studied by Lanteri-Turrini [{\bf LT}], who showed that $M$
is the Segre scroll over $\BP^1$ in this case; namely $S\cong\BP^1$,
$\SE\cong\SO(1)^{\oplus n}$, $M\cong\BP^1\times\BP^{n-1}$ and $M\subset\BP^N$
is the Segre embedding.
In this paper we are interested mainly in the case $s=2$.

This problem was studied by Ottaviani [{\bf Ot}] and
Beltrametti-Schneider-Sommese [{\bf BSS}], [{\bf BS}] in case $n=3$ and by
Ionescu-Toma [{\bf IT}] for general $n$.
Their results are as follows.

\Th (cf. [{\bf Ot}]). {\sl Let $M\subset\BP^5$ be a three-dimensional scroll
over a surface $S$.
Then one of the following conditions is satisfied.\nl
{\rm (1) (Segre scroll)} $S\cong\BP^2$, $\SE\cong\SO(1)^{\oplus 2}$, $d=\dg
M=3$.\nl
{\rm (2) (Bordiga scroll)} $S\cong\BP^2$, $c_1(\SE)=\SO(4)$, $c_2(\SE)=10$ and
$d=6$.\nl
{\rm (3) (Palatini scroll)} $S$ is isomorphic to a cubic surface,
$c_1(\SE)=\SO_S(2)$, $c_2(\SE)=5$ and $d=7$.
\nl
{\rm (4) (K3 scroll)} $S$ is a K3-surface obtained as a linear section of the
Grassman variety parametrizing lines in $\BP^5$ embedded by Pl\"ucker, and
$\SE$ is the restriction of the tautological vector bundle.
$c_2(\SE)=5$ and $d=9$.

All these four cases actually occur.}

\Th (cf. [{\bf IT}]). {\sl Let $M\subset\BP^N$ be a scroll over a surface $S$
with $N=2n-1$, $n>3$ and let $\SE$ and $d=\dg M$ be as above.
Then, $M$ is one of the following types:\nl
{\rm{($S_n$)}} $S\cong\BP^2$, $\SE\cong\SO(1)^{\oplus(n-1)}$,
$M\cong\BP^2\times\BP^{n-2}$, $d=n(n-1)/2$.\nl
{\rm{($B_n$)}} $S\cong\BP^2$, $c_1(\SE)=\SO(n+1)$, $c_2(\SE)=(n+1)(n+2)/2$,
$d=n(n+1)/2$.\nl
{\rm{($M_n$)}} $S$ is a K3 surface, $c_1(\SE)^2=2n^2-4$, $c_2(\SE)=n^2-4$,
$d=n^2$.\nl
{\rm{($?_n$)}} $S$ is a surface of general type.}

The existence of scrolls of the above type ($S_n$) is classical and is due to
Segre.
The case ($B_n$) is shown to exist for every $n$ in [{\bf IT}].
The existence of the type ($M_n$) is proved by Mukai ([{\bf Mu}]).
On the other hand, no example of type ($?_n$) is found:  [{\bf IT}] suspect
rather that there is no such scroll.
Any way, such scrolls must satisfy several numerical relations among their
invariants.
Here I propose a conjecture concerning further classifications derived from
these relations, which is verified for $n\le 11$ by hand and for $n\le 1100$ by
a computer.
In particular, a scroll of the type ($?_n$) does not exist unless $n=6, 10, 11,
12, 16, 18, 20, 24, 30, \cdots$.
On the other hand, no example is known and the existence problem is unsettled
for these $n$ (and also for $n>1100$, of course).

Mathematical tools used here are almost the same as [{\bf IT}], but we review
them here for the convenience of the reader.
When I started this study, I was not aware of this paper [{\bf IT}].
I would like to express my hearty thanks to Professors Ottaviani who informed
me of the result in [{\bf IT}].
I thank Professor Mukai, who communicated to me the existence of scrolls of the
type ($M_n$).
I also thank Professors Lanteri and Schneider for their cooperations in e-mail
correspondence during the preparation of this paper.
\dnl
{\bf \S1. Computing Chern classes}

(1.1) Throughout this paper let $M$ be a scroll in $\BP^{2n-1}$ over $S$ as in
\S0.
Let $\SE$ be the vector bundle on $S$ of rank $r=n-s+1$, $s=\dm S$, such that
$M\cong\BP_S(\SE)$ and the tautological bundle $H(\SE)$ is the hyperplane
section bundle $\SO_M(1)$, which will be denoted simply by $H$ from now on.
We put $d=\dg M=H^n\{M\}$.

(1.2) \Fact {\sl Via the ring homomorphism $\pi^*: H^\cdot(S)\to H^\cdot(M)$,
the cohomology ring of $M$ becomes a free $H^\cdot(S)$-module generated by $1,
h, h^2, \cdots, h^{r-1}$, where $h=c_1(H)\in H^2(M)$.
Moreover $\sum_{i=0}^r (-h)^i e_{r-i}=0$ in $H^\cdot(M)$, where $
e_j=\pi^*c_j(\SE)$.}

This is well known for general vector bundle $\SE$ of rank $r$.
Thus, the ring structure of $H^\cdot(M)$ is determined by $H^\cdot(S)$ and
Chern classes of $\SE$.

(1.3) \Fact {\sl Put $s_i(\SE)=\pi_*h^{r-1+i}\in H^{2i}(S)$ and
$s(\SE)=\sum_{i=0}^\infty s_i(\SE)\in H^\cdot(S)$.
Then $s(\SE)c(\SE^\vee)=1$, where $c(\SE^\vee)$ is the total Chern class of the
dual bundle $\SE^\vee$ of $\SE$.}

This is also standard.
$s_i(\SE)$ is called the $i$-th Segre class of $\SE$ and $s(\SE)$ is called the
total Sege class of $\SE$.
Moreover the following formulas are well known:
$$\align
s_1(\SE)&=c_1(\SE),\\
s_2(\SE)&=c_1(\SE)^2-c_2(\SE), \\
&\cdots\\
\endalign$$

(1.4) \Cor {\sl For any vector bundle $E$ on $X$ with $\rk E=r$ and for any
line bundle $L$ on $X$ with $c_1(L)=\ell$, we have $s_i(E\otimes
L)=\sum_j\binom{r-1+i}{j}s_{i-j}(E)\ell^j$.}

Indeed, $s_i(E\otimes
L)=\pi_*((h+\ell)^{(r-1+i)})=\sum_j\binom{r-1+i}{j}\pi_*h^{r-1+i-j}\cdot
\ell^j$ for $\pi:\BP(E)\to X$.

(1.5) Let $\pi: M=\BP(\SE)\to S$ be the projection.
Let $\SA$ be the kernel of the natural surjection $\pi^*\SE\to \SO[H]$.
Then we have an exact sequence $0\to\SO\to\SE^\vee\otimes H\to\SA^\vee\otimes
H\to 0$.
This is identified with the relative Euler sequence, and so we have the exact
sequence $0\to\SA^\vee\otimes H\to \Theta_M\to\pi^*\Theta_S\to 0$,
where $\Theta_X$ denotes the tangent bundle of $X$.

{}From these exact sequences we obtain the following relation
$c(\Theta_M)=\pi^*c(\Theta_S)c(\pi^*\SE^\vee\otimes H)$
of total Chern classes.

(1.6) Let $\SN$ be the normal bundle of $M$ in $\BP^{2n-1}$ and let $\Theta$ be
the restriction of the tangent bundle of $\BP^{2n-1}$ to $M$.
Then we have $c(\Theta_M)c(\SN)=c(\Theta)=(1+h)^{2n}$, where $h=c_1(H)$.

Combining (1.5), we get
$c(\SN)=(1+h)^{2n}c(\pi^*\SE^{\vee}\otimes H)^{-1}\pi^*c(\Theta_S)^{-1}
=(1+h)^{2n}s(\pi^*\SE\otimes[-H])\pi^*s(\Omega_S)$.

(1.7) In the following computation, $\pi^*\alpha$ will be denoted simply by
$\alpha$ for $\alpha\in H^\cdot(S)$.

By (1.4), we have
$s(\SE_M\otimes[-H])=\sum_{i=0}^\infty(\sum_{j=0}^{r-1+i}\binom{r-1+i}{j}s_{i-j}(\SE)(-h)^j)$.
Therefore the component of $(1+h)^{2n}s(\SE_M\otimes[-H])$ of degree $2k$ is
$\sum_{i,j}\binom{2n}{k-i}\binom{r-1+i}{j}(-1)^j h^{k-i+j}s_{i-j}(\SE)
\allowmathbreak=\sum_\ell(\sum_j\binom{2n}{k-j-\ell}\binom{r-1+j+\ell}{j}(-1)^j)h^{k-\ell}s_\ell(\SE)
\allowmathbreak=\sum_\ell\binom{2n-r-\ell}{k-\ell}h^{k-\ell}s_\ell(\SE)$
by the following

(1.8) {\it Claim}. $\sum_{(a,b)\vert
a+b=c}\binom{m-1+a}{a}(-1)^a\binom{p}{b}=\binom{p-m}{c}$.

To see this, use the Taylor expansion
$(1+T)^{-m}=\sum_{a\ge 0}\binom{m-1+a}{a}(-T)^a$
and compute the coefficients of $T^c$ in $(1+T)^{-m}(1+T)^p=(1+T)^{p-m}$.

(1.9) Combining (1.6) and (1.7), we get a formula for $c(\SN)$.
On the other hand, we have $c_n(\SN)=0$ and $c_{n-1}(\SN)=dh^{n-1}$ since $M$
is of codimension $n-1$ in $\BP^N$.
This gives non-trivial relations among Chern classes of $\SE$ and $\Theta_S$.
In the next section we analyse them precisely in case $s=\dm S=2$ and $r=n-1$.
\dnl
{\bf \S2. Over a surface}
\nopagebreak\par\nopagebreak
{}From now on we assume $s=\dm S=2$ and set $e_j=c_j(\SE)$,
$\gamma_i=c_i(\Theta_S)$.

(2.1) From (1.6) and (1.7) we obtain\nl
$c_{n-1}(\SN)=\binom{n+1}{n-1}h^{n-1}+\binom{n}{n-2}h^{n-2}e_1+\binom{n-1}{n-3}h^{n-3}(e_1^2-e_2)
-\binom{n+1}{n-2}h^{n-2}\gamma_1-\binom{n}{n-3}h^{n-3}e_1\gamma_1
+\binom{n+1}{n-3}h^{n-3}(\gamma_1^2-\gamma_2)$.
\nl
This is $dh^{n-1}$, while $h^{n-1}-h^{n-2}e_1+h^{n-3}e_2=0$ is the unique
relation in $H^\cdot(M)$ (cf. (1.2)).
Hence, substituting $h^{n-1}=h^{n-2}e_1-h^{n-3}e_2$ and comparing the
coefficients of $h^{n-2}$ and $h^{n-3}$, we get the following relations in
$H^\cdot(S)$:
$$(n^2-d)e_1=(n+1)n(n-1)\gamma_1/6{\text {,\qquad and}}\tag{i}$$
$$\tfrac{(n-1)(n-2)}2e_1^2+(d-n^2+n-1)e_2-\tfrac{n(n-1)(n-2)}6e_1\gamma_1+\tfrac{(n+1)n(n-1)(n-2)}{24}(\gamma_1^2-\gamma_2)=0.\tag{ii}$$
Next from $c_n(\SN)=0$ we get
$$3e_1^2-2e_2-ne_1\gamma_1+\tfrac{(n+1)(n-1)}{6}(\gamma_1^2-\gamma_2)=0.\tag{iii}$$

(2.2) Eliminating $\gamma_2$ from (ii) and (iii), we get
$-3(n^2-4)e_1^2+6(2d-n^2-2)e_2+n(n-2)(n+2)e_1\gamma_1=0$.
Eliminating $\gamma_1$ further using (i) and noting $d=(e_1^2-e_2)\{S\}$, we
obtain
$$\{2(q+2)d-q(q+5)\}\{(q+2)(q-4)d-(q+2)^2e_2+(q-1)(q-4)\}=-q(q-1)(q-4)(q+5)\tag{*}$$
for $q=n^2$.

(2.3) Thus, $2(q+2)d-q(q+5)$ is one of the finitely many divisors of the right
hand side, so there are only finitely many numerical possibilities for $d$ and
$e_2$.
Note that they are positive integers since $\SE$ is ample.

(2.4) For each $(d, e_2)$ satisfying (*), we examine the relation (i).
By the result [{\bf IT}], it suffices to consider the case where the canonical
bundle of $S$ is ample, or equivalently, $d>q=n^2$.
Let $6(d-n^2)/(n+1)n(n-1)=a/b$ for coprime integers $a, b$.
Then $ae_1\sim -b\gamma_1$, so $e_1\sim bA$ and $K_S\sim-\gamma_1\sim aA$ for
some ample line bundle $A$ on $S$.
Therefore
$$e_1^2=d+e_2 {\text{\quad is divided by \quad}} b^2.$$

(2.5) In the above case, we have $\gamma_1^2=a^2(d+e_2)/b^2$ and
$e_1\gamma_1=-a(d+e_2)/b$.
Using (iii), we solve $\gamma_2$.
The result has to satisfy the Noether relation:
$$\gamma_1^2+\gamma_2\equiv 0 {\text{\ modulo\ }} 12.$$

(2.6) It is easy to produce a computer programm to enumerate pairs $(n, d)$
satisfying the numerical conditions (2.2), (2.4) and (2.5).
In view of the result of our experiment, we make the following

\Conj {\sl Any pair $(n,d)$ with $n>3$ as above is one of the following
types:\nl
{\rm(1)} $n\equiv 0,2,6,12{\text{\ or\ }}16 {\text{\rm{\ modulo\ }}}18$ and
$d=q(q+5)/6$ for $q=n^2$.
Moreover $e_2=(q-4)(q+3)/6$, $K_S\sim ne_1$, $K_S^2=q(q^2+2q-6)/3$,
$\gamma_2=e(S)=(q^3+8q^2+24q+36)/3$ and $\chi(\SO_S)=(q^3+5q^2+9q+18)/18$ for
the corresponding scrolls.\nl
{\rm(2)} $n=10$ and $d=595$.
Moreover $e_2=561$, $K_S\sim 3e_1$, $K_S^2=10404$, $e(S)=12648$ and
$\chi(\SO_S)=1921$.\nl
{\rm(3)} $n=11$ and $d=231$.
Moreover $e_2=221$, $2K_S\sim e_1$, $K_S^2=113$, $e(S)=283$ and
$\chi(\SO_S)=33$.}

(2.7) {\it Remarks.}\nl
(1) The above conjecture is verified for $n\le 11$ by my hand, and for $n\le
1100$ by a personal computer.
In particular, the case with smallest $n$ is $(n,d)=(6,246)$.
[{\bf IT}] apparently claims that such a case with $n=6$ can be ruled out by
``divisibility manipulations'', but I cannot see how this can be done.\nl
(2) The above condition for $n$ of the type (1) is equivalent to $q\equiv
0{\text{\ or\ }}4{\text{\ modulo\ }}18$.\nl
(3) It is perhaps a delicate problem whether scrolls of the type ($?_n$) exist
or not for a pair $(n,d)$ in (2.6).
I find no example at present.
To settle the problem, we need some more geometric observations.
I feel that the case (2.6.3) might be of particular interest, since this is the
unique case with odd $n$ and moreover the invariants are relatively small.\nl
(4) The sectional genus $g=g(M,\SO(1))$ can be computed by using the relation
$2g-2=(K_S+e_1)e_1\{S\}$, and $g$ is bounded by the Castelnuovo inequality.
But it turns out that no pair $(n,d)$ in (2.6) is thus ruled out.

\enddocument